\documentclass[letterpaper, 10 pt, conference]{ieeeconf}  

\IEEEoverridecommandlockouts                              
\usepackage{graphics} 
\usepackage{epsfig} 
\usepackage{times} 
\usepackage{amsmath}

\newcommand\+{\mkern2mu}
\usepackage{amssymb}
\usepackage{physics}
\usepackage{mathtools}

\usepackage{algorithm}
\usepackage{algpseudocode} 

\usepackage{multicol}

\usepackage{subcaption}
\usepackage{caption}

\usepackage{xcolor}

\newcommand*{\QEDB}{\hfill\ensuremath{\blacksquare}}%

\newcommand*{\vsepfbox}[1]{%
  \begingroup
    \sbox0{\fbox{#1}}%
    \setlength{\fboxrule}{0pt}%
    \mbox{\kern-\fboxsep\fbox{\unhbox0}\kern-\fboxsep}%
  \endgroup
}

\title{\LARGE \bf
DN-ADMM: Distributed Newton ADMM for Multi-agent Optimization
}

\author{Yichuan Li$^{1}$, Nikolaos M. Freris$^{2}$, Petros Voulgaris$^{3}$ and Du\v{s}an Stipanovi\'{c}$^{1}$
\thanks{*This work was supported by the Ministry of Science and Technology of China under grant 2019YFB2102200, the Anhui Dept. of Science and Technology under grant 201903a05020049, the Tencent Holdings Ltd. under grant FR202003. }
\thanks{$^{1}$Coordinated Science Laboratory, University of Illinois at Urbana-Champaign, IL 61820, USA. {\tt\small yli129@illinois.edu}, {\tt\small dusan@illinois.edu}. }  %
\thanks{$^{2}$ School of Computer Science, University of Science and Technology of China, Hefei, Anhui, 230027, China. {\tt\small nfr@ustc.edu.cn.}}
\thanks{$^{3}$Department of Mechanical Engineering, University of Nevada, Reno, NV 89557, USA. {\tt\small pvoulgaris@unr.edu}.}
}

\begin{document}

\setlength{\belowdisplayskip}{2pt} 
\setlength{\abovedisplayskip}{2pt} 
\setlength{\textfloatsep}{0pt}

\maketitle
\thispagestyle{empty}
\pagestyle{empty}

\begin{abstract}
In a multi-agent network, we consider the problem of minimizing an objective function that is expressed as the sum of private convex and smooth functions, and a (possibly) non-differentiable convex regularizer.  We propose a novel distributed second-order method based on the framework of Alternating Direction Method of Multipliers (ADMM), by invoking approximate Newton iterations to the primal update corresponding to the differentiable part. In order to achieve a distributed implementation, the total Hessian matrix is split into a diagonal component (locally computable) and an off-diagonal component (that requires communication between neighboring agents). Subsequently, the Hessian inverse is approximated by a truncation of the Taylor expansion to $K$ terms: this amounts to fully distributed updates entailing $K$ distributed communication rounds. We establish global linear convergence to the primal-dual optimal solution under the assumption that the private functions are strongly convex and have Lipschitz continuous gradient. Numerical experiments demonstrate the merits of the approach comparatively with state-of-the-art methods.
\end{abstract}

\section{INTRODUCTION}
Decentralized consensus optimization problems consider an objective function consisting of a sum of cost functions which are only available at the corresponding agent of the underlying network. In many practical scenarios, a convex nonsmooth regularization function is augmented to the objective: typical examples include the $\ell_1-$norm and the nuclear norm that are known to promote sparsity in the solution structure \cite{c22}. Formally, let $\hat{x}\in \mathbb{R}^d$ denote the decision variable and consider a network containing $n$ agents. Each agent $i$ has access to the local private cost function $f^i:\mathbb{R}^d\to \mathbb{R}$ and exchanges information with neighbors so as to minimize the global objective function: 
\begin{gather}
    \hat{x}_\star = \underset{\hat{x}\in \mathbb{R}^d}{\text{argmin}} \,\bigg\{\sum_{i=1}^n f^i(\hat{x})+g(\hat{x}) \bigg\}.\label{1}
\end{gather} 
The global decision variable $\hat{x}$ is common to all agents but the local cost function $f^i(\cdot)$ is only available to agent $i$. Although it is possible to aggregate information from the network and proceed to obtain an approximate solution in a centralized fashion, it is desirable to solve problem (\ref{1}) through distributed exchanges. This is especially vital in view of the Big Data involved pertaining to large population count $n$, problem dimension $d$, and size of private data used in defining $f^i(\cdot)$, that are typically shared via limited bandwidth communication subject to the requirement for high responsiveness and real-time decision-making in Cyber-Physical Systems (CPS). In the distributed consensus optimization framework, every agent $i$ keeps a local copy of the decision variable $x^i$ and communicates with their neighbors to \emph{cooperatively} solve (\ref{1}). The setting has found multiple applications in distributed control \cite{c1}, power systems \cite{c2}, machine learning \cite{c3}, sensor networks \cite{c4}, and controller design \cite{c5}.


Proximal gradient \cite{c7} is a first-order method which generalizes the (projected) gradient method using the proximal operator to accommodate nonsmoothness. However, a distributed implementation is not possible when the nonsmooth term is not separable over the entries of the decision variable. Primal-dual methods for problem (\ref{1}) introduce local copies of the decision vector $x$ at each agent (i.e., agent $i$ maintains $x^i$) and introduce equality constraints between an agent's copy and its neighbors', which enforces consensus in the network \cite{c23}. The Alternating Direction Method of Multipliers (ADMM) \cite{c8} splits the primal variable over the smooth and the nonsmooth part of the objective function and sequentially obtains iterates by minimizing the augmented Lagrangian function. Another closely related method is the Method of Multipliers (MM) \cite{c24} which performs joint minimization over primal variables. However, the presence of nonsmooth part renders such approach cumbersome, in fact inapplicable in our setting. By completion of squares and the use of the proximal operator, the authors in \cite{c11} derived the {\it Proximal Augmented Lagrangian} which is continuously differentiable and eliminates the primal variable corresponding to the nonsmooth part. At each iteration of ADMM (and MM), it requires to solve an optimization sub-problem to obtain the next update, thus maybe especially costly in terms of computations. To avoid exact minimization steps required by the aforementioned primal dual algorithms, several approximation schemes have been proposed in the literature \cite{c12}-\cite{c14}. Decentralized Linearized ADMM (DLM) approximates the convex objective, at each iteration, with a first-order model and yields per-iteration costs comparable with gradient methods but suffers from slow convergence speed. Newton's method uses second-order information and features high convergence rates along with insensitivity to the condition number of the objective. However, Newton's iterations require solving linear systems as well as performing {\it backtracking} to ensure global convergence; this makes a distributed implementation impractical. 

We propose a Distributed Newton method in the framework of ADMM (DN-ADMM) and we establish global linear convergence without the need for backtracking. Compared to existing algorithms such as [14,15], we support nonsmooth regularizers using second order information without computing the generalized gradient of the proximal mapping. Numerical simulations demonstrate the advantage of our algorithm over existing methods. 

{\it Notation}: Vectors $x\in \mathbb{R}^d$ are represented as column vectors, denoted with lower case letters  and matrices $A\in \mathbb{R}^{n\times m}$ with capital letters. We use superscript to denote the vector component $x^i$ and subscript to denote the variable at the $t$-th iteration $x_t$. For a matrix $A$, we denote $A^{ij}$ as the $ij$-th entry of $A$, 
and when a norm subscript is not provided, $||x||$ and $||A||$ denote the Euclidean norm of a vector and the corresponding induced norm of a matrix, respectively. The norm of a vector $x$ with respect to a positive definite matrix $P$ is denoted as $||x||_P = (x^\top  Px)^{1/2}$ and the set $\{1,2,...,n\}$ is abbreviated  as $[n]$.

\section{PRELIMINARIES}
\subsection{Problem formulation}
Consider an undirected connected graph $\mathcal{G}=(\mathcal{V}, \mathcal{E}, \mathcal{A})$ of order $n$ with vertex set $\mathcal{V}=\{v^1,...,v^n\}$, edge set $\mathcal{E}\subseteq \mathcal{V}\times \mathcal{V}$, where each edge is denoted as $e^{ij}=(v^i,v^j)$, and the signed incidence matrix $\mathcal{A}\in \mathbb{R}^{n\times m}$, where $m=|\mathcal{E}|$ is the number of edges. We denote the neighborhood of vertex $i$ as $\mathcal{N}^i=\{v^j\in \mathcal{V}:(v^i,v^j)\in \mathcal{E}\}$. Problem (\ref{1}) can be reformulated to the standard consensus form for distributed optimization as follows. Let $x^i\in \mathbb{R}^d$ be the local copy of the decision variable held by the agent $i$ and introduce an extra variable $\hat{z}\in \mathbb{R}^d$:
\begin{gather}
    \underset{x^i,\,\hat{z}\in\mathbb{R}^d}{\text{minimize}}\,\,\left\{\sum_{i=1}^n f^i(x^i)+g(\hat{z})\right\} \nonumber\\
    \text{s.t.} \,\,x^i=x^j\,\,\quad\text{for all}\,\,i,\,\text{and}\,j\in \mathcal{N}^i,   \label{2}\\
    x^l = \hat{z} \,\,\quad\text{for some}\,\,l\in [n].\nonumber 
\end{gather}
If the underlying network graph $\mathcal{G}$ is connected, then the constraints in (\ref{2}) enforce network-wide consensus and (\ref{2}) is equivalent to (\ref{1}). The purpose of adding an extra variable $\hat{z}$ serves to separate the smooth and the nonsmooth part of the objective function and thus casts the problem to the ADMM framework. Note that in (\ref{2}), we only enforce equality constraint between $x^l$ and $\hat{z}$ for some agent $l$ in the network to minimize computation burden. Other frameworks such as enforcing the equality constraint for all agents have been considered as well and our approach is applicable therein with few modifications. We may express the above optimization problem compactly by aggregating $F(x)=\sum_{i=1}^n f^i(x^i)$, where we stack local copies into a vector $x=[(x^1)^\top  ,...,(x^n)^\top  ]^\top   \in \mathbb{R}^{nd}$, and further define $G(z)=g(\hat{z})$ where $z=[0,...,0,\hat{z}^\top  ]^\top  \in \mathbb{R}^{(m+1)d}$. Therefore, (2) becomes:  
\begin{gather}
\label{3}
    \underset{x,\,z}{\text{minimize}}\,\,F(x)+G(z) \\
    \text{s.t.} \,\, B^\top  x=z, \nonumber 
\end{gather} 
where we use Kronecker product to define $B=[A,E^l]\in \nolinebreak \mathbb{R}^{nd\times (m+1)d}$, $A=\mathcal{A}\otimes I_d$, and $E^l=e^l\otimes I_d$ is obtained by using the coordinate selection vector $e^l\in \mathbb{R}^n$, which is a zero vector except the $l$-th entry being equal to 1. We proceed to state our first assumption on the objective function. 

{\it Assumption 1}. The smooth local cost functions $\{f^i(\cdot)\}_{i=1}^n$ are twice differentiable with bounded Hessian as follows: 
\begin{gather}
    m_fI\preceq\nabla^2 f^i(x^i) \preceq M_fI, \quad \forall \,\,i\in [n].\label{4}
\end{gather}
The lower bound above implies local cost functions are strongly convex with parameter $m_f>0$ while the upper bound implies Lipschitz continuity of the gradient $\nabla f^i(\cdot)$ with constant $M_f>0$. Moreover, since $F(x)=\sum_{i=1}^n f^i(x^i)$, $\nabla^2 F(x)$ is a block-diagonal matrix with the $i$-th block being $\nabla^2 f^i(x^i)$. Consequently, the same bounds as in (\ref{4}) apply to the Hessian of the total cost function as well, i.e., $m_fI\preceq \nabla^2 F(x)\preceq M_fI$.
\subsection{Proximal operators}
Proximal operators can be seen as the generalization of the projection mapping for nonsmooth functions. For a closed, proper, and convex function $g:\mathbb{R}^d\to \mathbb{R}$, we define its proximal operator with parameter $\mu>0$ as:
\begin{gather}
\label{5a}
        \textbf{prox}_{\mu g}(v) = \underset{z}{\text{argmin}}\,\, \bigg\{g(z)+\frac{1}{2\mu}||v-z||^2\bigg\},
\intertext{and the associated value function, also known as the Moreau envelope,}
\label{5b}
M_{\mu g}(v) = g(\textbf{prox}_{\mu g}(v))+\frac{1}{2\mu}||v-\textbf{prox}_{\mu g}(v)||^2.
\end{gather}
The proximal mapping (\ref{5a}) is Lipschitz continuous with constant 1, differentiable almost everywhere, and single-valued \cite{c16}. 
For example, when $g(\cdot)$ is the indicator function of some convex closed set $\mathcal{X}$, its proximal operator reduces to the projection mapping of $\mathcal{X}$, i.e., 
$
    \Pi_{\mathcal{X}}(v) = \underset{z\in\mathcal{X}}{\text{argmin}} \norm{v-z}^2 
$.
For $g(z)=\sum \abs{z^i}$, the proximal operator has a closed form expression known as the soft-thresholding, i.e., $\textbf{prox}_{\mu g}(v)=\text{sign}(v)\,\text{max}(\abs{v}-\mu,0)$. More information on proximal operators and proximal algorithms can be found in \cite{c7}.
\subsection{Alternating Direction Method of Multipliers}
The augmented Lagrangian associated with (\ref{3}) with parameter $\mu>0$ is obtained by augmenting a quadratic penalty term of the constraint to the Lagrangian:
$
    \mathcal{L}_{\mu}(x,z;y) = F(x)+G(z)+ y^\top  (B^\top  x-z)
    +(1/2\mu)\norm{B^\top  x-z}^2 ,
$
where $y\in \mathbb{R}^{(m+1)d}$ is the dual variable associated with constraints. ADMM takes advantage of the fact that primal variables are separated among smooth and nonsmooth parts of the cost function and sequentially updates them by minimizing $\mathcal{L}_\mu (x,z;y)$ as follows:
\begin{subequations}
\label{8}
\begin{gather}
    x_{t+1} = \underset{x}{\text{argmin}}\,\, \mathcal{L}_{\mu}(x,z_t;y_t) ,\label{8a}\\
    z_{t+1} = \underset{z}{\text{argmin}}\,\, \mathcal{L}_{\mu}(x_{t+1}, z;y_t) ,\label{8b}\\ 
    y_{t+1} = y_t+\tfrac{1}{\mu} \nabla_y \mathcal{L}(x_{t+1},z_{t+1};y_t). \label{8c}
\end{gather}
\end{subequations}
In the following section, we propose a variant of ADMM method with an additional primal update and perform inexact minimization for step (\ref{8a}). We extend the result in \cite{c14} and \cite{c20} to include a nonsmooth regularizer by building a quadratic model at each iteration with Hessian inverse approximated by truncated Taylor expansion, while invoking proximal operators to update the $z$-variable.
\section{Second-order approximation and distributed implementation}
\subsection{Second order approximation }
Newton's method can be interpreted as minimizing a quadratic model incurred by the Hessian of the objective function, at each iteration. However, backtracking is required to ensure global convergence and damped stepsizes are often used \cite{c21}. A quadratic model for (\ref{8a}) can be constructed as $\hat{\mathcal{L}}_{\mu}(x,z_t;y_t)\approx \mathcal{L}_\mu(x_t,z_t;y_t)+\nabla_x \mathcal{L}_\mu(x_t,z_t;y_t)^\top  (x-x_t) 
    +\frac{1}{2}(x-x_t)^\top  \nabla_{xx}^2\mathcal{L}_\mu(x_t,z_t;y_t)(x-x_t)+\tfrac{\epsilon}{2}\norm{x-x_t}^2$.
The last term $\tfrac{\epsilon}{2}||x-x_t||^2$ aims to keep the next iterate close to the current one and therefore renders a numerically robust algorithm. We obtain the primal update for $x_{t+1}$ as:
\begin{gather}
    x_{t+1} = x_t - H_t^{-1}[\nabla F(x_t)+\tfrac{1}{\mu}BB^\top  x_t-\tfrac{1}{\mu}Bz_t+By_t]. \nonumber 
\end{gather}
The Hessian $H_t=\nabla_{xx}\hat{\mathcal{L}}_{\mu}$ is given by:
\begin{gather}
    H_t = \nabla^2 F(x_t)+\tfrac{1}{\mu}BB^\top  +\epsilon I_{nd}, \label{17}
\end{gather}
where $I$ is the identity matrix of appropriate dimension. \\
\indent Recall that $e^l$ is the coordinate selection vector and since $\nabla^2 F(x_t)$ is block-diagonal and $BB^\top  =(L+e^l(e^l)^\top  )\otimes I_d$ (where $L$ is the graph Laplacian), the Hessian $H_t$ is \emph{distributedly computable} by using local information and communication with neighbors. \\
\indent However, the Hessian inverse $H_t^{-1}$ is not necessarily distributedly computable. Motivated by \cite{c14} and \cite{c17}, we employ a matrix splitting technique so that the Hessian inverse is suitable for Taylor expansion and, therefore, we can truncate the series to obtain a finite-term approximation of the Hessian inverse.
We first decompose the Hessian as $H_t=D_t-N$ where $D_t$ is block diagonal and $N$ encodes off-diagonal entries. In specific, $D_t \succ0$ is equal to:
\begin{gather}
\label{Dt}
    D_t = \nabla^2 F(x_t)+\tfrac{1}{\mu}(e^l(e^l)^\top  )\otimes I_d+\tfrac{2}{\mu}L_{\text{dia}}\otimes I_d+\epsilon I_{nd}, 
\end{gather}
where $e^l\in \mathbb{R}^n$ and $N\succeq0$:
\begin{gather}
    N = \tfrac{1}{\mu}(L_{\text{dia}}-L_{\text{off}})\otimes I_{d},
\end{gather}
where $L_{\text{dia}}$ and $L_{\text{off}}$ are the diagonal and off-diagonal components of the graph Laplacian matrix $L$. With this matrix splitting, we can expand the Hessian inverse as $H_t^{-1} = D_t^{-\frac{1}{2}}(I-D_t^{-\frac{1}{2}}ND_t^{-\frac{1}{2}})^{-1}D_t^{-\frac{1}{2}}$, and use Taylor expansion to express $(I-D_t^{-\frac{1}{2}}ND_t^{-\frac{1}{2}})^{-1}$ as $\sum_{i=0}^{\infty} (D_t^{-\frac{1}{2}}ND_t^{-\frac{1}{2}})^i$. We truncate the series to finite many terms and show that the error incurred is bounded. To ensure validity of this expansion, we need to show that the spectral radius of $D_t^{-\frac{1}{2}}ND_t^{-\frac{1}{2}}$ is strictly less than $1$ which we prove next. \\
\indent {\it Lemma 1}. The spectral radius of $D_t^{-\frac{1}{2}}ND_t^{-\frac{1}{2}}$ is strictly less than $1$. \\
\indent {\it Proof :} The following proof is inspired by \cite{c14}. We first define $C = \tfrac{2}{\mu}L_{\text{dia}}\otimes I\succ0$. 
We bound the spectral radius by considering upper bounds of each term as $\norm{D_t^{-\frac{1}{2}}ND_t^{-\frac{1}{2}}} \leq \norm{D_t^{-\frac{1}{2}}C^{\frac{1}{2}}}^2\norm{C^{-\frac{1}{2}}NC^{-\frac{1}{2}}} $. For the first term, note that $D_t^{-\frac{1}{2}}C^{\frac{1}{2}}$ is a block-diagonal matrix with $i$-th block specified as:
$$
    (D_t^{-\frac{1}{2}}C^{\frac{1}{2}})^{ii} = 
    (\nabla^2 f^i(x^i_t)+\tfrac{1}{\mu}((e^l)^i +2L_{\text{dia}}^{ii}+\mu\epsilon)I)^{-\frac{1}{2}}\nonumber\\$$
\begin{gather}
    \times(\tfrac{2}{\mu}L_{\text{dia}}^{ii})^{\frac{1}{2}} 
     = \left( \frac{\mu(\nabla^2 f^i(x^i_t)+\epsilon I)+(e^l)^{i}I }{2L_{\text{dia}}^{ii}}+I\right)^{-\frac{1}{2}},
\end{gather}
From {\it Assumption 1} and the fact that $1\leq (L_{\text{dia}})^{ii}\leq n-1$ for any $i$, we can bound $(D_t^{-\frac{1}{2}}C^{\frac{1}{2}})^{ii}$ as:
\begin{gather}
   (D_t^{-\frac{1}{2}}C^{\frac{1}{2}})^{ii}  \preceq \bigg( \frac{2(n-1)}{\mu(m_f+\epsilon)+2(n-1)}\bigg)^{\frac{1}{2}} I. \label{12}
\end{gather}
Since $C^{-\frac{1}{2}}NC^{-\frac{1}{2}}=C^{-\frac{1}{2}}(NC^{-1})C^\frac{1}{2}$, we can bound the eigenvalues of $C^{-\frac{1}{2}}NC^{-\frac{1}{2}}$ by investigating the eigenvalues of $NC^{-1} =\frac{1}{2}(L_{\text{dia}}-L_{\text{off}})L_{\text{dia}}^{-1}\otimes I$. Note that $L_{\text{dia}}$ is diagonal with entries being equal to the number of neighbors of the corresponding agent, while $L_{\text{off}}$ is symmetric with $ij$-th entry being equal to $-1$ if and only if $(v^i,v^j) \in \mathcal{E}$ and $0$ otherwise. Therefore,
\begin{gather}
    NC^{-1} = \frac{1}{2} 
    \label{22}
    \begin{bmatrix}    
    1  & \frac{-L_{\text{off}}^{12}}{L_{\text{dia}}^{22}}  &\dots & \frac{-L_{\text{off}}^{1n}}{L_{\text{dia}}^{nn}} \\
    \frac{-L_{\text{off}}^{21}}{L_{\text{dia}}^{11}} & 1 &\dots &   \frac{-L_{\text{off}}^{2n}}{L_{\text{dia}}^{nn}} \\
    \vdots & \vdots &  \ddots & \vdots \\ 
    \frac{-L_{\text{off}}^{n1}}{L_{\text{dia}}^{11}} & \frac{-L_{\text{off}}^{n2}}{L_{\text{dia}}^{22}} &\dots & 1\\
    \end{bmatrix}\otimes I.
\end{gather}
The sum of non-diagonal components of column $i$ is 
$
    \frac{1}{2}\sum_{j=1,j\neq i}^n  \frac{-L_{\text{off}}^{ji}}{L_{\text{dia}}^{ii}} =\frac{1}{2},
$
where we have used the property of graph Laplacian $\sum_{j} L^{ji}=0$ for any $i$. Since each diagonal entry of the above matrix equals $\frac{1}{2}$, by applying Gershgorin circle theorem, we get that all eigenvalues of the matrix shown in (\ref{22}) lie in the circle centered at $\frac{1}{2}$ with radius $\frac{1}{2}$. Furthermore, for matrices $A\in \mathbb{R}^{n\times n}$ and $B\in\mathbb{R}^{m\times m}$, the eigenvalues of $A\otimes B$ are given by $\lambda_i(A)\lambda_j(B)$ for $1\leq i\leq n$ and $1\leq j\leq m$ where $\lambda_i$ stands for eigenvalues. We conclude that 
$
    0\leq \lambda_i(NC^{-1}) \leq 1.  
$
 Therefore, combining with (\ref{12}), we obtain an upper bound for the spectral radius of $D_{t}^{-\frac{1}{2}}ND_t^{-\frac{1}{2}}$ as follows:
\begin{gather}
\label{lemma1}
    \sigma_{max} (D_{t}^{-\frac{1}{2}}ND_t^{-\frac{1}{2}}) \leq \frac{2(n-1)}{\mu(m_f+\epsilon)+2(n-1)}<1.
\end{gather}
\QEDB
\vspace{-5mm}
\subsection{Distributed implementation}
We proceed to develop DN-ADMM with Hessian inverse approximated by truncating the Taylor series and show its distributed implementation in Algorithm 1. Define the $K$-th order Hessian inverse approximation as:
\begin{gather}
    \widehat{H}_{t}^{-1}(K) = D_t^{-\frac{1}{2}}\sum_{i=0}^K(D_t^{-\frac{1}{2}}ND_t^{-\frac{1}{2}})^iD_t^{-\frac{1}{2}}. \label{25}
\end{gather}
Since only the last block of $z_t$ is nonzero, step (\ref{8b}) of ADMM is equivalent to:
$
    \hat{z}_{t+1} = \underset{\hat{z}}{\text{argmin}}\{g(\hat{z})+(y_t^{m+1})^\top (x^l_{t+1}-\hat{z})+\tfrac{1}{2\mu}||x^l_{t+1}-\hat{z}||^2\}
$.
By completion of squares, this update reduces to the proximal mapping of $g(\cdot)$ evaluated at $(x^l_{t+1}+\mu y_t^{m+1})$ with parameter $\mu$, i.e.,
\begin{gather}
\hat{z}_{t+1}= \textbf{prox}_{\mu g}(x^l_{t+1}+\mu y_t^{m+1}).
\end{gather}
Note that we have used the updated value $x_{t+1}^l$ and the current dual $y_t^{m+1}$ to obtain $\hat{z}_{t+1}$. In DN-ADMM, we perform an extra proximal mapping step before updating $x_{t+1}$ to mimic the behavior of the method of multipliers. The updates of order $K$ (number of terms kept in the Taylor expansion) of DN-ADMM are given below:
\begin{subequations}
\begin{align}
\hat{\theta}_{t+1} &= \textbf{prox}_{\mu g}(x_t^l+\mu y_t^{m+1}) ,\label{27a}\\
    x_{t+1} &= x_t-\hat{H}_t^{-1}(K)[\nabla F(x_t)+\tfrac{1}{\mu}B(B^\top x_t-\theta_{t+1})\nonumber\\&\qquad+By_t]  ,\label{27b}\\
    \hat{z}_{t+1} &=\textbf{prox}_{\mu g}(x_{t+1}^l+\mu y_t^{m+1}) , \label{27c}\\
    y_{t+1} &= y_t+\tfrac{1}{\mu}(B^\top x_{t+1}-z_{t+1}) ,\label{27d}
\end{align}\label{27} 
\end{subequations}\noindent where $\theta_{t+1}\in \mathbb{R}^{(m+1)d}$ is a zero vector except the last blocking being equal to $\hat{\theta}_{t+1}$, similarly for $\hat{z}_{t+1}$ and $z_{t+1}$. In contrary to \cite{c11} and \cite{c15}, proximal operators do not explicitly appear in the $x$-update and therefore the approximated Hessian $\hat{H}_t$ does not include the generalized gradient of proximal mapping. This results in a more efficient implementation, that avoids the need for computing (potentially cumbersome) generalized Hessians pertinent to the non-differentiable part. Moreover, $K$-th order approximation of the Hessian inverse requires $K$ rounds of communication among agents. Specifically, define $u_t(K-1) = \hat{H}^{-1}_t(K-1)\nabla_x \mathcal{L}_\mu(x_t,\theta_{t+1};y_t)$ so that the $x$-updates is simply $x_{t+1}=x_t-u_t(K-1)$. From (\ref{25}), we express the next order update direction $u_t(K)$ as:
\begin{align}
    &u_t(K)=D_t^{-\frac{1}{2}}\sum_{i=0}^{K}(D_t^{-\frac{1}{2}}ND_t^{-\frac{1}{2}})^iD_t^{-\frac{1}{2}}\nabla_x\mathcal{L}_\mu  \nonumber\\
    &=D_t^{-1}\nabla_x \mathcal{L}_\mu+D_t^{-1}ND_t^{-\frac{1}{2}}\sum_{i=0}^{K-1}(D_t^{-\frac{1}{2}}ND_t^{-\frac{1}{2}})^iD_t^{-\frac{1}{2}}\nabla_x \mathcal{L}_\mu \nonumber\\
    \label{28}
    &=D_t^{-1}\nabla_x\mathcal{L}_\mu+D_t^{-1}N u_t(K-1), 
\end{align} 
\noindent where we have suppressed the argument of $\nabla_x \mathcal{L}_\mu(x_t,z_t;y_t)$ for brevity. Note that since $D_t$ is block diagonal and $N$ respects the network structure, we can compute $u_t(K)$ by just one round of distributed communication once $u_t(K-1)$ is obtained. Therefore, to obtain the $K$-th order approximation of the Hessian inverse, every agent initializes with $u_t(0)=D_t^{-1}\nabla_x\mathcal{L}_\mu$, then communicates for $K$ rounds with its neighbors. Moreover, since $\theta_{t+1},z_{t+1}$ has first $m$ blocks of entries being zero, only the $(m+1)$-th block, i.e., the subvector $\hat{\theta}_{t+1},\hat{z}_{t+1}\in \mathbb{R}^d$, needs to be computed, only using  information held by the $l$-th agent.\\
\begin{algorithm}[t]
    \caption{DN-ADMM($K$) updates at the agent $i$} 
    Zero initialization $x_0^i$, $\phi_0^i$ and $y^{m+1}_0$. Hyperparameters $\mu$, $\epsilon$.
    \begin{algorithmic}[1]
    \For{$t=0,1,2,\ldots$}
    \State Compute the local block $D_t^{ii}$ using (\ref{Dt}).
    \State Compute $h_t^i=\nabla f^i(x_t^i)+\tfrac{1}{\mu}\sum_{j\in \mathcal{N}^i} (x_t^i-x_t^j) + \phi_t^i$.
            \If{$i=l$}
                \State Compute $\hat{\theta}_{t+1}=\textbf{prox}_{\mu g}(x_t^i+\mu y^{m+1}_t)$.
                \State Update $h_t^i\leftarrow h_t^i +\tfrac{1}{\mu}(x_t^i-\hat{\theta}_{t+1})+y^{m+1}_t$.
            \EndIf
    \State Form $u^i_t(0)=(D_t^{ii})^{-1}h_t^i$.
                    \For{$k=0,1,2,\hdots K-1$}
                    \State $u^i_t(k+1)=(D_t^{ii})^{-1}h_t^i+\tfrac{1}{\mu}\sum_{j\in\mathcal{N}^i}(u_t^i(k)+u_t^j(k))$.
                    \EndFor
    \State \textbf{Primal update:} $x_{t+1}^i=x_t^i-u_t^i(K)$.
    \State \textbf{Dual update:} $  \phi_{t+1}^i = \phi_t^i+\tfrac{1}{\mu}\sum_{j\in \mathcal{N}^i}(x_{t+1}^i-x_{t+1}^j)$.
    \If{$i=l$}
    \State Compute $\hat{z}_{t+1}=\textbf{prox}_{\mu g}(x^i_{t+1}+\mu y^{m+1}_t)$.
    \State Compute $y^{m+1}_{t+1} = y^{m+1}_t +\tfrac{1}{\mu}(x^i_{t+1}-\hat{z}_{t+1})$.
    \EndIf\\
    \textbf{Information exchange}\\
	Agent $i$ exchanges $x^i_t,x_{t+1}^i,u_t^i(k)$ with neighbors.
    \EndFor
\end{algorithmic}
\end{algorithm}
\indent Updates shown in (\ref{27}) follow directly from the framework of ADMM and are convenient for analysis purpose. We further present an efficient distributed implementation in Algorithm 1. From the definition of $B$, we note that in (\ref{27b}), the dual variable $y_t$ is used in the update of $x_{t+1}$ in the form of $By_t= A y^{[m]}_t+E^l y^{m+1}_t$, where $y^{[m]}_t\in \mathbb{R}^{md}$ denotes the concatenation of the first $m$ blocks of $y_t$. By defining $\phi_t=Ay_t^{[m]}\in \mathbb{R}^{nd}$, we eliminate the need to exchange dual variables during the primal update and obtain a distributed implementation as follows. We let each agent hold the corresponding pair $(x^i,\phi^i)$ while the $l$-th agent (which updates the nonsmooth variable) additionally holds $y^{m+1}$. At each iteration, agent $i$ begins by computing the block $D_t^{ii}$ using local information and proceeds to compute $h_t^i:=\nabla_x \mathcal{L}^i_\mu$ using values $x_t^j$ obtained from its neighbors (step 3). The $l$-th agent additionally evaluates the proximal mapping associated with the nonsmooth regularizer $g(\cdot)$ and updates $h_t^i$ accordingly in step 6. Steps 9-11 amount to $K$ rounds of communication among the network to compute $u_t(K)$ as in (\ref{28}). Agents proceed to perform primal and dual updates in steps 12 and 13 respectively while the $l$-th agent additionally performs proximal mapping with updated iterates and update the additional dual variable $y^{m+1}_{t+1}$ in steps 15 and 16 respectively. 
\section{CONVERGENCE ANALYSIS}
In this section, we prove that iterates generated by (\ref{27}) converge linearly to the optimal $x_\star,\,z_\star$ and $y_\star$. We first state the KKT conditions of the problem along with an additional assumption on cost functions needed for our analysis in the following. 

{\it KKT conditions} for problem (\ref{3}):
\begin{align*}
    \nabla F(x_\star)+By_\star &= 0 \tag{KKTa}\label{KKTa}\\
    \partial g(\hat{z}_\star) - y^{m+1}_\star &\ni 0 \tag{KKTb}\label{KKTb}\\
    B^\top  x_\star-z_\star &=0 \tag{KKTc} \label{KKTc}
\end{align*}

{\it Assumption 2}. The function $g(z)$ is proper, closed, and convex, i.e., $\forall w_1\in \partial g(z_1), w_2\in\partial g(z_2)$,
\begin{gather}
    (z_1-z_2)^\top  (w_1-w_2)\geq 0. 
\end{gather}
\indent {\it Lemma 2}. Consider the dual updates specified in (\ref{27d}), with zero initialization. Then $y_t$ is in the range of $B^\top  $ for all $t$. Moreover, there is a unique dual optimal $y_\star$ in the column space of $B^\top$.

{\it Proof} : See Appendix. 

The following lemma characterizes a property of the dual iterates in relation to the nonsmooth variable $\hat{z}_{t+1}$.\\
\indent {\it Lemma 3}. Consider the update $\hat{z}_{t+1}=\textbf{prox}_{\mu g}(x_{t+1}^l+\mu y_t^{m+1})$ in (\ref{27c}). It holds that: 
\begin{gather}
    y_{t+1}^{m+1}\in \partial g(\hat{z}_{t+1}). \label{inclusion}
\end{gather}
{\it Proof} : See Appendix. \\
\indent {\it Lemma 4}. Consider the updates of $x_{t+1}$ and $y_{t+1}$ specified in (\ref{27}), and recall the approximated Hessian (\ref{25}). For primal/dual optimal solutions $x_\star$ and $y_\star$, it holds that: 
\begin{gather}
\nabla F(x_{t+1})-\nabla F(x_\star)+B(y_{t+1}-y_\star)
+e_t \nonumber\\ +\epsilon(x_{t+1}-x_t) =0\,,  \label{lemma4}
\end{gather}
where $e_t$ is defined as:
\begin{gather}
        e_t = \nabla F(x_t)+\nabla^2 F(x_t)(x_{t+1}-x_t)-\nabla F(x_{t+1}) \nonumber\\
    +(\hat{H}_t-H_t)(x_{t+1}-x_t)+\tfrac{1}{\mu}B(z_{t+1}-\theta_{t+1}). \label{error}
\end{gather}
{\it Proof} : See Appendix. \\
\indent {\it Lemma 5}. Consider the error term defined in (\ref{error}), it is bounded above as:
\begin{gather}
\label{lemma5}
\norm{e_t} \leq \gamma \norm{x_{t+1}-x_t},
\end{gather}
where $\gamma$ is defined as: 
\begin{gather}
\scalebox{0.9}{$
   \gamma = 2M_f+\tfrac{1}{\mu}
    +(M_f+\epsilon+\tfrac{1+2(n-1)}{\mu})\left(\frac{2(n-1)}{\mu(m_f+\epsilon)+2(n-1)}\right)^{K+1}$}. \label{gamma}
\end{gather}
{\it Proof:} See Appendix. \\
\indent Consider the following vector matrix defined in terms of primal/dual variables and algorithm hyperparameters:
\begin{gather}
v = \begin{bmatrix}x \\y \end{bmatrix}\quad \mathcal{H} = \begin{bmatrix} \epsilon I&0\\0 &\mu I \end{bmatrix}. \label{34}
\end{gather}
Note that $v$ is the concatenation of $(x,y)$ and $\mathcal{H}$ is a diagonal matrix with blocks being corresponding scaling parameters. By defining $v_\star$ as the concatenation of the optimal variables, i.e., $v_\star = [x_\star^\top  ,y_\star^\top  ]^\top  $, we prove the convergence of the iterates to the optimality by showing that the Lyapunov function $\norm{v_t-v_\star}^2$ converges to zero. \\
\indent {\it Theorem 1}. Consider the iterates generated by DN-ADMM in (\ref{27}). Denote the smallest positive eigenvalue of $B^\top B$ as $\lambda_\mathrm{min}$. Let $\beta,\eta >1,\zeta\in (\tfrac{m_f+M_f}{2 m_fM_f},\tfrac{\epsilon}{\gamma^2})$ be arbitrary constants and choose $\epsilon >\left[2M_f+\tfrac{m_f+2M_f(n-1)}{\mu m_f}\right]^2\tfrac{m_f+M_f}{2m_fM_f}$. Recall the definition of $v$ and $\mathcal{H}$ in (\ref{34}). If Assumptions 1, 2 hold, then the sequence $\{\norm{v_{t}-v_\star}^2_{\mathcal{H}}\}_{t\in\mathbb{N}}$ converges linearly:
\begin{gather}
    \norm{v_{t+1}-v_\star}^2_{\mathcal{H}} \leq \frac{1}{1+\delta}\norm{v_t-v_\star}^2_{\mathcal{H}}, \label{theorem1}
\end{gather}
where  
\begin{gather}
    \delta = \min\bigg\{\tfrac{2m_fM_f}{\epsilon(m_f+M_f)}-\tfrac{1}{\epsilon\zeta},\tfrac{(\epsilon-\zeta\gamma^2)\lambda_{\mathrm{min}}(\beta-1)(\eta-1)}{\mu\beta[\epsilon^2(\eta-1)+\eta\gamma^2(\beta-1)]},\nonumber\\
    \tfrac{2\lambda_{\mathrm{min}}}{(m_f+M_f)\mu\beta\eta}\bigg \}.\label{choose}
\end{gather}
\indent {\it Proof} : Since $F(x)$ is strongly convex with parameter $m_f$ and the gradient $\nabla F(x)$ is Lipschitz continuous with parameter $M_f$, it holds that \cite{cbook}: 
\begin{gather*}\scalebox{0.9}{$
    \frac{m_fM_f}{m_f+M_f}\norm{x_{t+1}-x_\star}^2 
    +\frac{1}{m_f+M_f}\norm{\nabla F(x_{t+1})-\nabla F(x_\star)}^2  \leq $}
\end{gather*}
\begin{gather}
    (x_{t+1}-x_\star)^\top  (\nabla F(x_{t+1})-\nabla F(x_\star)) . \label{30}
\end{gather}
We substitute the expression for $\nabla F(x_{t+1})-\nabla F(x_\star)$ from {\it Lemma 4} and denote the left-hand-side of (\ref{30}) by \textbf{LHS}; this can be upper bounded by:
\begin{gather}
    \textbf{LHS} \leq -(x_{t+1}-x_\star)^\top  B(y_{t+1}-y_\star)-(x_{t+1}-x_\star)^\top  e_t\nonumber\\
    -\epsilon (x_{t+1}-x_\star)^\top  (x_{t+1}-x_t) . \nonumber
\end{gather}
From the dual update (\ref{27}) and KKTc, we have
$
     -(x_{t+1}-x_\star)^\top  B = -\mu(y_{t+1}-y_t)^\top  -(z_{t+1}-z_\star)^\top             
$.
Further substituting this expression for $-(x_{t+1}-x_\star)^\top  B$ into above inequality gives:
\begin{gather}
    \textbf{LHS} \leq -\mu(y_{t+1}-y_t)^\top  (y_{t+1}-y_\star)-\nonumber\\
    (z_{t+1}-z_\star)^\top  (y_{t+1}-y_\star)-(x_{t+1}-x_\star)^\top  e_t-\nonumber\\
    \epsilon(x_{t+1}-x_\star)^\top  (x_{t+1}-x_t) \label{39}
\end{gather}
Since $(z_{t+1}-z_\star)^\top (y_{t+1}-y_\star)=(\hat{z}_{t+1}-\hat{z}_\star)^\top(y_{t+1}^{m+1}-y_t^{m+1})$, using {\it Lemma 3}, KKTb, and {\it Assumption 2}, we have 
$
    (\hat{z}_{t+1}-\hat{z}_\star)^\top  (y_{t+1}^{m+1}-y_\star^{m+1}) \in (\hat{z}_{t+1}-\hat{z}_\star)^\top   (\partial g(\hat{z}_{t+1})-\partial g(\hat{z}_\star ))\geq 0.
$
Therefore, \textbf{LHS} can be further bounded as:
\begin{gather}
    \textbf{LHS} \leq -\mu(y_{t+1}-y_t)^\top  (y_{t+1}-y_\star)-(x_{t+1}-x_\star)^\top  e_t\nonumber\\
    -\epsilon(x_{t+1}-x_\star)^\top  (x_{t+1}-x_t). \nonumber
\end{gather}
We use the identity $-2(a-b)^\top  (a-c)=\norm{b-c}^2-\norm{a-b}^2-\norm{a-c}^2$ with $a=y_{t+1}, b=y_t$ and $c=y_\star$ and similarly, with $(x_{t+1}, x_t,x_\star)$. Multiplying the above inequality by $2$ on both sides yields:
\begin{gather}
    \norm{x_{t+1}-x_\star}^2_{\frac{2m_fM_f}{m_f+M_f}I}+\norm{\nabla F(x_{t+1})-\nabla F(x_\star)}^2_{\frac{2}{m_f+M_f}I}\leq \nonumber\\
    \mu\big(\norm{y_t-y_\star}^2-\norm{y_{t+1}-y_t}^2-\norm{y_{t+1}-y_\star}\big)\nonumber\\
    +\epsilon\big(\norm{x_t-x_\star}^2
    -\norm{x_{t+1}-x_t}^2-\norm{x_{t+1}-x_\star}^2\big)\nonumber\\
    -2(x_{t+1}-x_\star)^\top  e_t. \label{41}
\end{gather}
Since $2(x_{t+1}-x_\star)^\top  e_t\geq -\frac{1}{\zeta}\norm{x_{t+1}-x_\star}^2-\zeta\norm{e_t}^2\+\+\forall\,\zeta>0$ and using the definition (\ref{34}), we rewrite (\ref{41}) as:
\begin{gather}
    \norm{x_{t+1}-x_\star}^2_{(\frac{2m_fM_f}{m_f+M_f}-\frac{1}{\zeta})I}
    +\norm{\nabla F(x_{t+1})-\nabla F(x_\star)}^2_{\frac{2}{m_f+M_f}I}\nonumber \\
    +\mu\norm{y_{t+1}-y_t}^2+\epsilon\norm{x_{t+1}-x_t}^2-\zeta\norm{e_t}^2 \nonumber\\
    \leq \norm{v_t-v_\star}^2_\mathcal{H}-\norm{v_{t+1}-v_\star}^2_\mathcal{H}.   \label{42}
\end{gather}
From the dual update in (\ref{27d}), the first $m$ block of the dual variable is updated as $y_{t+1}^{[m]}=y^{[m]}_t+\tfrac{1}{\mu}A^\top x_{t+1}$. After rearranging and combining with the first $m$ block of KKTc, we have $ \norm{y_{t+1}-y_t}^2\geq \norm{y_{t+1}^{[m]}-y_t^{[m]}}^2=\tfrac{1}{\mu^2}\norm{x_{t+1}-x_\star}^2_{L\otimes I_d}$, where we have used the fact that $L\otimes I_d=AA^\top$. To prove linear convergence, we need to show that for some $\delta>0,\,\norm{v_{t+1}-v_\star}^2_{\mathcal{H}}\leq \frac{1}{1+\delta}\norm{v_t-v_\star}^2_{\mathcal{H}}$. Therefore, along with the lower bound for $\norm{y_{t+1}-y_t}^2$ substituted in, we show that the \textbf{LHS} of (\ref{42}) is lower bounded by $\delta\norm{v_{t+1}-v_\star}^2_{\mathcal{H}}$, i.e., 
\begin{gather}
   \norm{x_{t+1}-x_\star}^2_{(\frac{2m_fM_f}{m_f+M_f}-\frac{1}{\zeta})I+\tfrac{1}{\mu}L\otimes I_d}
   +\epsilon\norm{x_{t+1}-x_t}^2\nonumber \\ +\norm{\nabla F(x_{t+1})-\nabla F(x_\star)}^2_{\frac{2}{m_f+M_f}I}-\zeta\norm{e_t}^2
   \nonumber\\
    \geq \delta\norm{v_{t+1}-v_\star}^2_{\mathcal{H}}\label{43}
\end{gather}
We proceed to establish such a bound by using the components of the \textbf{LHS} of (\ref{43}). From (\ref{lemma4}), we have:
\begin{gather}
    B(y_{t+1}-y_\star)=-\big[\nabla F(x_{t+1})-\nabla F(x_\star)+\epsilon(x_{t+1}-x_t)\nonumber\\
    +e_t  \big] \label{44}
\end{gather}
Furthermore, it is easy to verify that when $a=b+c$ implies that $\forall\,\beta>1$, it holds that $\norm{a}^2\leq \frac{\beta}{\beta-1}\norm{b}^2+\beta\norm{c}^2$. Applying this formula two times with respective constants $\beta,\eta>1$, we obtain an upper bound on the \textbf{LHS} of (\ref{44}) as:
\begin{gather*}
    \norm{B(y_{t+1}-y_\star)}^2\leq \tfrac{\beta\epsilon^2}{\beta-1}\norm{x_{t+1}-x_t}^2+\tfrac{\beta\eta}{(\eta-1)}\norm{e_t}^2\nonumber \\
    +\beta\eta\norm{\nabla F(x_{t+1})-\nabla F(x_\star)}^2,
\end{gather*}
From {\it Lemma 2}, it follows that dual iterates are orthogonal to the kernel of $B$ and therefore, $\norm{B(y_{t+1}-y_\star)}^2\geq \lambda_{\mathrm{min}}\norm{y_{t+1}-y_\star}^2$, where $\lambda_{\mathrm{min}}$ is the smallest positive eigenvalue of $B^\top  B$. Therefore, we have the following inequality: 
\begin{gather}
    \norm{y_{t+1}-y_\star}^2\leq \tfrac{\beta\epsilon^2}{\lambda_{\mathrm{min}}(\beta-1)}\norm{x_{t+1}-x_t}^2+\tfrac{\beta\eta}{\lambda_{\mathrm{min}}(\eta-1)}\norm{e_t}^2\nonumber\\
    +\tfrac{\beta\eta}{\lambda_{\mathrm{min}}}\norm{\nabla F(x_{t+1})-\nabla F(x_\star)}^2.\label{33}
\end{gather}
Combining upper bounds for $\norm{y_{t+1}-y_\star}^2,\norm{e_t}$ in (\ref{33}) and (\ref{lemma5}) respectively, from the definition of $\norm{v_{t+1}-v_\star}^2_\mathcal{H}=\mu\norm{y_{t+1}-y_\star}^2+\epsilon\norm{x_{t+1}-x_\star}^2$, we obtain:
\begin{gather*}
    \norm{v_{t+1}-v_\star}^2_{\mathcal{H}}\leq (\tfrac{\mu\beta\epsilon^2}{\lambda_{\mathrm{min}}(\beta-1)}+\tfrac{\mu\beta\eta\gamma^2}{\lambda_{\mathrm{min}}(\eta-1)})\norm{x_{t+1}-x_t}^2\\
    +\tfrac{\mu\beta\eta}{\lambda_{\mathrm{min}}}\norm{\nabla F(x_{t+1})-\nabla F(x_\star)}^2+\epsilon \norm{x_{t+1}-x_\star}^2.
\end{gather*}
Therefore, to satisfy (\ref{43}), it is sufficient to show for some $\delta>0$,
\begin{gather}
    \norm{x_{t+1}-x_\star}^2_{(\frac{2m_fM_f}{m_f+M_f}-\frac{1}{\zeta})I+\tfrac{1}{\mu}L\otimes I_d}
   +\epsilon\norm{x_{t+1}-x_t}^2 \nonumber\\ +\norm{\nabla F(x_{t+1})-\nabla F(x_\star)}^2_{\frac{2}{m_f+M_f}I} \geq
    \nonumber\\
     (\tfrac{\delta\mu\beta\epsilon^2}{\lambda_{\mathrm{min}}(\beta-1)}+\tfrac{\delta\mu\beta\eta\gamma^2}{\lambda_{\mathrm{min}}(\eta-1)}+\zeta\gamma^2)\norm{x_{t+1}-x_t}^2\nonumber\\
    \tfrac{\delta\mu\beta\eta}{\lambda_{\mathrm{min}}}\norm{\nabla F(x_{t+1})-\nabla F(x_\star)}^2+\delta\epsilon \norm{x_{t+1}-x_\star}^2.
  \label{51}
\end{gather}
Inequality (\ref{51}) holds if $\delta$ is chosen to satisfy:
\begin{subequations}
\begin{gather}
\left(\tfrac{2m_fM_f}{m_f+M_f}-\tfrac{1}{\zeta}-\delta\epsilon \right)I+\tfrac{1}{\mu}L\otimes I_d\succeq 0\\
\epsilon -\tfrac{\delta\mu\beta\epsilon^2}{\lambda_{\mathrm{min}}(\beta-1)}-\tfrac{\delta\mu\beta\eta\gamma^2}{\lambda_{\mathrm{min}}(\eta-1)}-\zeta\gamma^2 \geq 0 \\
\tfrac{2}{m_f+M_f}-\tfrac{\delta\mu\beta\eta}{\lambda_{\mathrm{min}}}\geq 0
\end{gather} \label{52}
\end{subequations}\noindent Inequalities (\ref{52}) are satisfied if $\delta$ chosen as in (\ref{choose}) and $\zeta,\epsilon$ as in the statement (which guarantees that $\delta>0$ in view of the first and the second quantity in (\ref{choose})). Therefore, inequality (\ref{43}) is satisfied and we conclude (\ref{theorem1}).\QEDB



\begin{figure}[t]
  \centering
  \begin{subfigure}[b]{0.49\linewidth}
    \includegraphics[width=\textwidth]{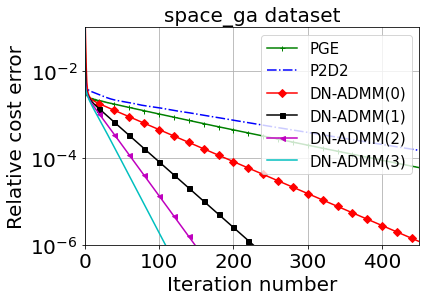}
    \caption{$n = 20,\,d=6$.}
  \end{subfigure}
  \begin{subfigure}[b]{0.49\linewidth}
    \includegraphics[width=\textwidth]{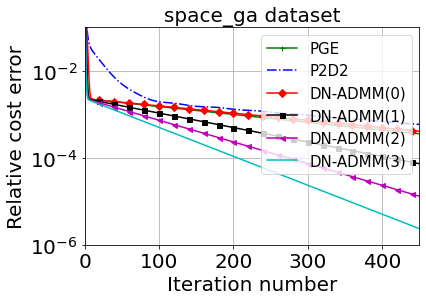}
    \caption{$n=40,\,d=6$.}
  \end{subfigure}
  \caption{space\_ga dataset.}
  \label{fig1}
  \begin{subfigure}[b]{0.49\linewidth}
    \includegraphics[width=\textwidth]{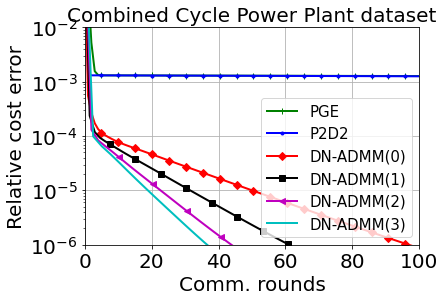}
    \caption{$n=30,\,d=4$.}
  \end{subfigure}
  \begin{subfigure}[b]{0.49\linewidth}
    \includegraphics[width=\textwidth]{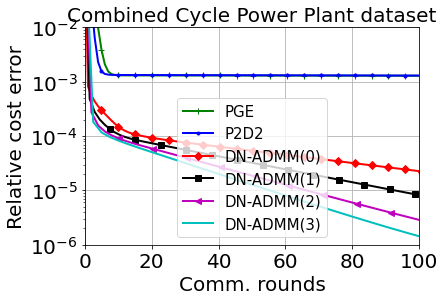}
    \caption{$n=50,\,d=4$.}
  \end{subfigure}
  \caption{Combined Cycle Power Plant dataset.}
  \label{fig2}
\end{figure}
\section{NUMERICAL SIMULATIONS}
We apply DN-ADMM($K$) to a distributed LASSO problem where each agent holds a local quadratic cost function with a common $\gamma$-weighted $\ell_1-$norm regularizer,
\begin{gather*}
    \underset{x\in\mathbb{R}^d}{\text{minimize}}\,\,l(x)=\bigg\{ \sum_{i=1}^n \frac{1}{2}\norm{A_ix-b_i}^2+\gamma\norm{x}_1\bigg\}.
\end{gather*}
We generate a connected random binomial graph of $n$ agents where an edge is drawn i.i.d Bernoulli($p$) with $p=0.2$. Two real datasets are considered from the LIBSVM\footnote{https://www.csie.ntu.edu.tw/~cjlin/libsvm/} and UCI Machine Learning Repository\footnote{https://archive.ics.uci.edu/ml/index.php}: the space\_ga dataset and Combined Cycle Power Plant dataset (CCPP). The space\_ga dataset contains observations on U.S county votes on the 1980 presidential election with 6 features and the CCPP dataset predicts the net hourly electrical energy output of a power plant with 4 features. We took 3,080 instances from the space\_ga dataset and 9,000 instances from the CCPP dataset, and evenly distributed it among $n$ agents. In Figs. \ref{fig1} and \ref{fig2}, we plot the averaged relative cost $\frac{\tfrac{1}{n}\sum_{i=1}^n l(x_t^i)-l(x_\star)}{\tfrac{1}{n}\sum_{i=1}^n l(x_0^i)-l(x_\star)}$ versus the number of iterations and the number of communication rounds per node. We compare the performance of between DN-ADMM($K$) with state-of-the-art first-order methods: P2D2 \cite{cp2d2} and PG-EXTRA \cite{c19} with the mixing matrix generated using the Metropolis rule and the Laplacian based constant edge weight matrix, respectively. We do not compare against second order methods since existing ones do not support nonsmooth regularizers. For all simulations, $\gamma=0.002$, stepsizes for PG-EXTRA and P2D2 are hand tuned to achieve fastest convergence. We observe significant speedup of DN-ADMM($K$) over first order methods from Fig.\ref{fig1} and Fig. \ref{fig2}. Moreover, we observe that by communicating with neighbors during the stage of Hessian approximation, DN-ADMM($K$) is able to achieve significantly faster convergence speed. Note that, in all cases of DN-ADMM($K$), we only communicate vectors exactly as first order methods while in DN-ADMM(0), \emph{no inner loop communication} is required to compute the update direction once the gradient of the augmented Lagrangian is obtained.


\section{CONCLUSIONS} 

We have proposed a Newton variant of ADMM for distributed consensus optimization, that applies Newton steps in the primal update rule corresponding to the differentiable component of the objective. The method relies on the proximal mapping of the nonsmooth part of the cost function and approximates the Hessian inverse via truncation of Taylor's series. Our analysis has established global linear convergence for the primal-dual iterates under the standard strong convexity and Lipschitz gradient assumptions without backtracking.

\appendix
\noindent \textbf{A. Proof of Lemma 2}\\
 Recall the dual update in (\ref{27}) and the fact that only the $(m+1)$-th block of $z_{t+1}$, i.e., $\hat{z}_{t+1}$, is not identically zero. Therefore, it suffices to show that $[\mathbf{0}^\top,\hat{z}_{t+1}^\top]^\top $ is in the range of $B^\top  $. Note that since the underlying graph is connected, $A^\top  b=\mathbf{0}$ for any $b\in \mathbb{R}^{nd}$ such that all $n$ blocks of sub-vectors of dimension $d$ are equal. Construct such a vector $b\in \mathbb{R}^{nd}$ as:
$
    b = \begin{bmatrix}
    \hat{z}_{t+1}^\top & \cdots & \hat{z}_{t+1}^\top
    \end{bmatrix}^\top \label{d1}
$. From the definition of $B$:
\begin{gather}
    B^\top  b = \begin{bmatrix}A^\top  b \\
    (E^l)^\top  b\end{bmatrix} = \begin{bmatrix}\mathbf{0} \\ \hat{z}_{t+1} \end{bmatrix}=z_{t+1}, \label{d2}
\end{gather}
which shows that $z_{t+1}$ is in the range of $B^\top$. Therefore, with zero initialization, dual iterates generated by (\ref{27}) stay in the range of $B^\top$. Next we proceed to show that there exists a unique dual optimal $y_\star$ in the range of $B^\top$. For any dual optimal $y_o$ that satisfy the KKTa, $\nabla F(x_\star)+By_o=0$, its projection to the column space of $B^\top$, denoted as $y_\star$, also satisfy KKTa. This is because their difference is in the kernel of $B$, i.e., $B(y_o-y_\star)=0$. We show uniqueness of $y_\star$ by contradiction. Assume there are two dual optimal solutions, $y_1=B^\top  b_1$ and $y_2=B^\top  b_2$ and $y_1\neq y_2$. From KKTa (uniqueness of $x_\star$ is guaranteed by strong convexity of $F(\cdot)$), we obtain:
\begin{gather*}
    \nabla F(x_\star)+ BB^\top  b_1 = 0 \\
    \nabla F(x_\star)+ BB^\top  b_2 = 0.
\end{gather*}
After taking the difference, we obtain $BB^\top  (b_1-b_2)=0$. Since $BB^\top  =L\otimes I_d+E^l(E^l)^\top  \succ 0$, we conclude $b_1=b_2$, contradiction.\QEDB

\noindent \textbf{B. Proof of Lemma 3}\\
Rearranging the dual updates in (\ref{27}) and focusing on the last block, we get:
\begin{gather}
    \hat{z}_{t+1} = x_{t+1}^l-\mu(y_{t+1}^{m+1}-y_t^m) \label{35}
\end{gather}
Moreover, since $\hat{z}_{t+1} =\textbf{prox}_{\mu g}(x^l_{t+1}+\mu y_t^{m+1})= \text{argmin}_v \{g(v)+\frac{1}{2\mu}\norm{x^l_{t+1}+\mu y^{m+1}_t-v}^2\}$, the optimality condition is
$
    0 \in \partial\+ g(\hat{z}_{t+1})+\tfrac{1}{\mu}(\hat{z}_{t+1}-x_{t+1}^l-\mu y_t^{m+1}) 
$.
Substituting the expression for $z_{t+1}$ from (\ref{35}) into above and canceling identical terms gives (\ref{inclusion}). \QEDB

\noindent \textbf{C. Proof of Lemma 4}
Consider the primal updates for $x_{t+1}$ in (\ref{27}). After re-arranging, we obtain:
\begin{gather}
\nabla F(x_t)+\tfrac{1}{\mu}B(B^\top  x_t-\theta_{t+1})+By_t+\hat{H}_t(x_{t+1}-x_t)=0. \nonumber
\end{gather}
Recalling the definition of the Hessian $H_t=\nabla^2 F(x_t)+\tfrac{1}{\mu}BB^\top  +\epsilon I$ in (\ref{17}), add and subtract $H_t(x_{t+1}-x_t)$ from the above gives,
\begin{gather*}
    \nabla F(x_t)+\nabla^2 F(x_t)(x_{t+1}-x_t)+By_t+\tfrac{1}{\mu}BB^\top  x_{t+1} \nonumber\\
    +(\hat{H}_t-H_t)(x_{t+1}-x_t)-\tfrac{1}{\mu}B\theta_{t+1}
    +\epsilon(x_{t+1}-x_t)=0.
\end{gather*}
With the definition of error term (\ref{error}), we rewrite the above as:
\begin{gather}
        e_t+\nabla F(x_{t+1})+By_t+\tfrac{1}{\mu}BB^\top  x_{t+1}-\tfrac{1}{\mu}Bz_{t+1}\nonumber \\
    +\epsilon(x_{t+1}-x_t) = 0 .
\end{gather}
Substituting dual updates into above and subtract KKTa, we obtain claimed (\ref{lemma4}). 
\QEDB 

\noindent \textbf{D. Proof of Lemma 5}\\
Recall the definition of error term in (\ref{error}) and apply triangle inequality to obtain:
\begin{gather}
    \norm{e_t}\leq \norm{\nabla F(x_t)+\nabla^2F(x_t)(x_{t+1}-x_t)-\nabla F(x_{t+1})}\nonumber\\
    +\norm{(\hat{H}_t-H_t)(x_{t+1}-x_t)}+\tfrac{1}{\mu}\norm{B(z_{t+1}-\theta_{t+1})}.
 \label{61}
\end{gather}
We proceed to bound each component separately. From Assumption 1 and Lipschitz continuity of the gradient:
\begin{gather}
    \norm{\nabla F(x_t)+\nabla^2F(x_t)(x_{t+1}-x_t)-\nabla F(x_{t+1})} \leq \nonumber \\
    \label{comp1}
    2M_f \norm{x_{t+1}-x_t}.
\end{gather}
For the second term in (\ref{61}), we first rewrite $\hat{H}_t-~H_t=\hat{H}_t^{\frac{1}{2}}(I-\hat{H}_t^{-\frac{1}{2}}H_t\hat{H}_t^{-\frac{1}{2}})\hat{H}_t^{\frac{1}{2}}$. Therefore, we have:
\begin{gather}
    \norm{\hat{H}_t-H_t}\leq \norm{\hat{H}_t^{\frac{1}{2}}}^2\norm{I-\hat{H}_t^{-\frac{1}{2}}H_t\hat{H}_t^{-\frac{1}{2}}}.\label{matrix1}
\end{gather}
From properties of similar matrices, it follows that:
\begin{gather}
\label{matrix2}
\norm{I-\hat{H}_t^{-\frac{1}{2}}H_t\hat{H}_t^{-\frac{1}{2}}}=\norm{I-H_t\hat{H}_t^{-1}}.
\end{gather}
Recall the matrix splitting $H_t=D_t-N$ and the definition in (\ref{25}) (recall also that $\hat{H}_t$ is abbreviation for $\hat{H}_t(K)$),
\begin{align}
    H_t\hat{H}_t^{-1} &= (D_t-N)(D_t^{-\frac{1}{2}}\sum_{i=0}^K(D_t^{-\frac{1}{2}}ND_t^{-\frac{1}{2}})^i D_t^{-\frac{1}{2}}) \nonumber \\
    &=(D_t-N)(\sum_{i=0}^KD_t^{-1}(ND_t^{-1})^i). \nonumber
\end{align}
Therefore, we have:
\begin{align}
    I-H_t\hat{H}_t^{-1} &= I-\sum_{i=0}^K(ND_t^{-1})^i+\sum_{i=0}^K(ND_t)^{i+1} \nonumber \\
    &=(ND_t^{-1})^{K+1}.
\end{align}
From (\ref{lemma1}) we can bound this as:
\begin{gather}
    \norm{I-H_t\hat{H}_t^{-1}}=\norm{ND_t^{-1}}^{K+1}=\norm{D_t^{-\frac{1}{2}}ND_t^{-\frac{1}{2}}}^{K+1}\nonumber\\
    \label{matrix3}
    \leq \left(\frac{2(n-1)}{\mu(m_f+\epsilon)+2(n-1)}\right)^{K+1}.
\end{gather}
Combining (\ref{matrix1}), (\ref{matrix2}) and (\ref{matrix3}), we have:
\begin{gather}
\label{63}
    \norm{\hat{H}_t-H_t}\leq \left(\frac{2(n-1)}{\mu(m_f+\epsilon)+2(n-1)}\right)^{K+1}\norm{\hat{H}_t^\frac{1}{2}}^2. 
\end{gather}
By definition, we know:
\begin{align}
\hat{H}_t^{-1}&=D_t^{-\frac{1}{2}}\sum_{i=0}^K(D_t^{-\frac{1}{2}}ND_t^{-\frac{1}{2}})^iD_t^{-\frac{1}{2}} \nonumber \\
&=D_t^{-1}+D_t^{-\frac{1}{2}}\sum_{i=1}^K(D_t^{-\frac{1}{2}}ND_t^{-\frac{1}{2}})^iD_t^{-\frac{1}{2}} . \label{def1} 
\end{align}
Since the second term in (\ref{def1}) is positive semidefinite, we obtain $\norm{\hat{H}^{-1}_t}\geq \norm{D_t^{-1}}$ and using a term-by-term upper bound on the eigenvalues of $D_t$ in (\ref{Dt}) we obtain: 
$$
    \norm{\hat{H_t}^{-1}}\geq \norm{D_t^{-1}}\geq \frac{\mu}{\mu(M_f+\epsilon)+\big[1+2(n-1)\big]} \nonumber \\
    $$
\begin{gather}
    \norm{\hat{H}_t}\leq M_f+\epsilon+\tfrac{1+2(n-1)}{\mu}. \label{65}
\end{gather}
Combining (\ref{63}) and (\ref{65}), we bound the second term of (\ref{61}) as:
\begin{gather}
    \norm{(\hat{H}_t-H_t)(x_{t+1}-x_t)}\leq \nonumber \\
    (M_f+\epsilon+\tfrac{1+2(n-1)}{\mu})\norm{x_{t+1}-x_t}\bigg(\tfrac{2(n-1)}{\mu(m_f+\epsilon)+2(n-1)}\bigg)^{K+1}. \label{66}
\end{gather}
The third term of (\ref{61}) can be upper bounded by considering the fact that $\tfrac{1}{\mu}\norm{B(z_{t+1}-\theta_{t+1})}=\tfrac{1}{\mu}\norm{\hat{z}_{t+1}-\hat{\theta}_{t+1}}=\tfrac{1}{\mu}\norm{\textbf{prox}_{\mu g}(x^l_{t+1}+\mu y^{m+1}_t)-\textbf{prox}_{\mu g}(x^l_t+\mu y^{m+1}_t)}$. Since the proximal operator is nonexpanisve, it follows that $\tfrac{1}{\mu}\norm{\hat{z}_{t+1}-\hat{\theta}_{t+1}}\leq \tfrac{1}{\mu}\norm{x_{t+1}-x_t}$. Combining these upper bounds , we arrive at the claim in (\ref{lemma5}). \QEDB\\

\end{document}